\documentclass[letterpaper, 10 pt, conference]{ieeeconf}  
\IEEEoverridecommandlockouts                              
\title{\LARGE \bf
	A one-size-fits-all artificial pancreas for people with type 1 diabetes based on physiological insight and feedback control
}

\author{Tobias K. S. Ritschel, Asbj{\o}rn Thode Reenberg, Emilie B. Lindkvist, Christian Laugesen, Jannet Svensson, \\ Ajenthen G. Ranjan, Kirsten N{\o}rgaard, Bernd Dammann, John Bagterp J{\o}rgensen
	\thanks{This work was partially funded by the IFD Grand Solution project ADAPT-T2D (9068-00056B). A. T. Reenberg, T. K. S. Ritschel, B. Dammann, and  J. B. J{\o}rgensen are with the Department of Applied Mathematics and Computer Science, Technical University of Denmark, DK-2800 Kgs. Lyngby, Denmark.
	E. B. Lindkvist, C. Laugesen, J. Svensson, A. G. Ranjan, and K. N{\o}rgaard are with Steno Diabetes Center Copenhagen, Clinical Research, DK-2730 Herlev, Denmark.
	Corresponding author: J. B. J{\o}rgensen (E-mail: {\tt\small jbjo@dtu.dk}).}
}

\usepackage{amsmath}
\usepackage{graphicx}
\usepackage{cite}
\usepackage{xcolor}

\begin{document}
	\maketitle
\thispagestyle{empty}
\pagestyle{empty}

	\begin{abstract}
	We propose a model-free artificial pancreas (AP) for people with type 1 diabetes. The algorithmic parameters are tuned to a virtual population of 1,000,000 individuals, and the AP repeatedly estimates the basal and bolus insulin requirements necessary for maintaining normal blood glucose levels. Therefore, the AP can be used without healthcare personnel or engineers customizing the algorithm to each user. The estimates are based on bodyweight, measurements from a continuous glucose monitor (CGM), and estimates of the meal carbohydrate contents.
	In a virtual clinical trial with all 1,000,000 individuals (i.e., a Monte Carlo closed-loop simulation), the AP achieves a mean time in range of more than 87\% and almost 89\% of the participants satisfy several glycemic targets.
\end{abstract}

	\section{Introduction}
\label{sec:introduction}
Diabetes is a chronic disease where the pancreas produces insufficient amounts of insulin or the body is resistant to insulin. More than 10\% of the world's adult population suffers from this disease, and in 2021, USD 966 billion dollars were spent on diabetes (which corresponds to 9\% of the global health expenditure)~\cite{IDF:2021}.
Type 1 diabetes (T1D) accounts for between 5--10\% of all cases, and it is caused by autoimmune destruction of the pancreatic insulin-producing cells. Consequently, the pancreas
does not produce any insulin. Therefore, people with T1D require daily insulin treatment in order to prevent high blood glucose concentrations (referred to as hyperglycemia).
Long periods of hyperglycemia can cause damage to the nerves and eyes and lead to chronic kidney disease and cardiovascular disease. Additionally, incorrect insulin treatment can lead to low blood glucose concentrations (referred to as hypoglycemia). Severe hypoglycemia can lead to a number of acute complications, e.g., loss of consciousness and seizures.

People with T1D spend significant amounts of time on self-treatment. Therefore, there is considerable interest in developing automated insulin delivery systems which can assist them. Such systems are referred to as \emph{artificial pancreases} (APs)~\cite{Lal:etal:2019}, and they typically consist of 1)~a sensor, often a continuous glucose monitor (CGM), 2)~a control system, usually a control algorithm implemented on a smartphone or a dedicated device, and 3)~an actuator, e.g., an insulin pump.

Many control algorithms have been considered for this purpose. They can be divided into model-free and model-based algorithms. Model-based algorithms typically use model predictive control (MPC)~\cite{Boiroux1:etal:2018, Boiroux2:etal:2018, Chakrabarty:etal:2020, Messori:etal:2018} where a model is used to predict the body's response to, e.g., meal carbohydrates and insulin. MPC is a well-proven control methodology that has been applied to many different types of processes~\cite{Forbes:etal:2015}. However, it requires an accurate model. Automatic generation of such a model based on historical data (e.g., CGM measurements, administered insulin, and meal carbohydrates) is an ongoing field of research. Consequently, it remains difficult to make model-based APs widely available because a model must be developed for each individual.
In contrast, model-free controllers only rely on a few pieces of information about the body for which estimates are readily available, e.g., the bodyweight, the basal insulin requirement, the insulin-to-carb ratio (ICR), and the insulin sensitivity factor (ISF). These controllers are often based on heuristics~\cite{Capel:etal:2014}, fuzzy logic~\cite{Biester:etal:2019}, or proportional-integral-derivative (PID) control.
PID controllers have been successfully applied in many different industrial applications~\cite{Astrom:Murray:2008}, and several researchers have proposed APs based on concepts from PID control. Marchetti et al.~\cite{Marchetti2:etal:2008} proposed a PID controller which is switched off when a meal is announced and switched back on based on heuristical rules. Huyett et al.~\cite{Huyett:etal:2015} described a PID control algorithm for intraperitoneal insulin delivery (whereas most APs deliver insulin subcutaneously which results in a more delayed insulin effect). Finally, J{\o}rgensen et al.~\cite{Jorgensen:etal:2019} and Sejersen et al.~\cite{Sejersen:etal:2021} present APs that 1)~use concepts from PID control to estimate the basal rate and 2)~compute the insulin bolus as a linear function of the meal carbohydrate content using the ICR (which is assumed to be known). However, these AP algorithms have not been tested on large numbers of real or virtual people. Therefore, it is currently unknown whether they can be adopted without healthcare personnel or engineers tuning the algorithms specifically for each individual (or group of individuals).

In this work, we present a one-size-fits-all AP algorithm with a single set of parameters tuned to a population of 1,000,000 virtual individuals with T1D. It simultaneously estimates the basal insulin and the meal insulin bolus curve (as a function of the meal carbohydrate content normalized with body weight). Therefore, it is straightforward for a user to start using the system. Furthermore, we demonstrate that, for a given objective function, the optimal meal insulin bolus is a nonlinear function of the meal carbohydrate content, and we argue that it can be approximated well by a piecewise linear function.
We test the proposed AP using a virtual clinical trial with 1,000,000 participants over 52 weeks, i.e., we perform a Monte Carlo closed-loop simulation. All virtual participants meet the target on time in level 2 hypoglycemia, the mean time in range (TIR) is 87.2\%, and almost 89\% of the virtual participants meet all targets on TIR, time above range (TAR), time below range (TBR), average glucose, and glucose management indicator (GMI) described in~\cite{Holt:etal:2021b}.

The remainder of this paper is organized as follows. In Section~\ref{sec:analysis}, we analyze the optimal insulin bolus as a function of the meal carbohydrate content, and in Section~\ref{sec:algorithm}, we present the AP. We present the results of the virtual clinical trial in Section~\ref{sec:numerical:example}, and conclusions are given in Section~\ref{sec:conclusions}. 
	\section{Analysis}
\label{sec:analysis}
In this section, we present a dynamic optimization problem for determining the optimal meal insulin bolus as a function of the meal size. We only use this optimization problem to analyze the meal insulin bolus curve. It is not used in the AP algorithm presented in Section~\ref{sec:algorithm} because it would require a model of each person using the AP.

\subsection{The dynamic optimization problem}
\label{sec:analysis:dynamic:optimization:problem}
The dynamic optimization problem determining the optimal meal insulin bolus flow rate is in the form
\begin{subequations}
	\label{eq:bolus}
	\begin{align}
		\label{eq:bolus:objective:function}
		\min_{u_0} \quad
		\phi &= \int_{t_0}^{t_f} \rho(z(t))\,dt,
	\end{align}
	subject to
	\begin{align}
		\label{eq:bolus:initial:condition}
		x(t_0) &= x_0, \\
		\label{eq:bolus:dynamical:constraint}
		\dot x(t) &= f(x(t), u(t), d(t), \theta), & t &\in [t_0, t_f], \\
		\label{eq:bolus:outputs}
		z(t) &= g(x(t), \theta), & t &\in [t_0, t_f], \\
		\label{eq:bolus:zero:order:hold:parametrization:u}
		u(t) &= u_k, \quad t \in [t_k, t_{k+1}[, & k &= 0, \ldots, N-1, \\
		\label{eq:bolus:zero:order:hold:parametrization:d}
		d(t) &= d_k, \quad t \in [t_k, t_{k+1}[, & k &= 0, \ldots, N-1, \\
		\label{eq:bolus:bounds}
		u_{\min} & \leq u_0 \leq u_{\max}.
	\end{align}
\end{subequations}
The objective function in~\eqref{eq:bolus:objective:function} is the integral of the penalty function $\rho$ over the time horizon $[t_0, t_f] = [0, 12]$~h, and $N$ is the number of control intervals. Furthermore, $t$ is time, $x$ are the state variables, $u$ is a vector of manipulated inputs (i.e., the basal and bolus insulin flow rates), $d$ are disturbance variables (i.e., the meal carbohydrate flow rate), $z$ is the output (i.e., the CGM measurement), and $\theta$ are model parameters. The initial condition~\eqref{eq:bolus:initial:condition} is the steady state corresponding to a blood glucose concentration of~6~mmol/L, and the dynamical constraint~\eqref{eq:bolus:dynamical:constraint} is the model presented by Hovorka et al.~\cite{Hovorka:etal:2004} extended with a CGM model~\cite{Facchinetti:etal:2014}. Next, \eqref{eq:bolus:outputs} is an output equation, and~\eqref{eq:bolus:zero:order:hold:parametrization:u}--\eqref{eq:bolus:zero:order:hold:parametrization:d} are zero-order-hold~(ZOH) parametrizations of the manipulated inputs and the disturbance variables. Finally, $u_{\min}$  and $u_{\max}$ in~\eqref{eq:bolus:bounds} are lower and upper bounds on the manipulated inputs in the first control interval.

In the first control interval ($k = 0$), the person consumes a meal with a specified meal carbohydrate flow rate, $d_0$, and an insulin bolus is administered. For the remaining control intervals ($k > 0$), the disturbances and the bolus insulin flow rate are zero. The insulin basal rate is equal to its steady state value (corresponding to 6~mmol/L) in all control intervals, i.e., it is not a decision variable.
Finally, the lower bound on the insulin bolus flow rate in the first control interval ($k = 0$) is 0 and the upper bound is infinity.

The penalty function penalizes the deviation from the setpoint $\bar z = 6$~mmol/L and the violation of the soft lower bound $z_{\min} = 3.9$~mmol/L (see also Table~\ref{tab:glycemic:ranges}), i.e.,
\begin{align}
	\label{eq:bolus:penalty:function}
	\rho(z(t)) &= \bar \rho(z(t)) + \kappa \rho_{\min}(z(t)),
\end{align}
where
\begin{subequations}
	\begin{align}
		\label{eq:bolus:setpoint:penalty}
		\bar \rho  (z(t)) &= \frac{1}{2} (z(t) - \bar z)^2, \\
		\label{eq:bolus:lower:bound:penalty}
		\rho_{\min}(z(t)) &= \frac{1}{2} \max\{0, z_{\min} - z(t)\}^2.
	\end{align}
\end{subequations}
As the main priority is to avoid hypoglycemia, $\kappa = 10^6$.

\subsection{Optimal insulin bolus curves}
\label{sec:analysis:optimal:boluses}
\begin{figure*}
	\centering
	\includegraphics[width=\textwidth, trim=140pt 35pt 100pt 35pt, clip]{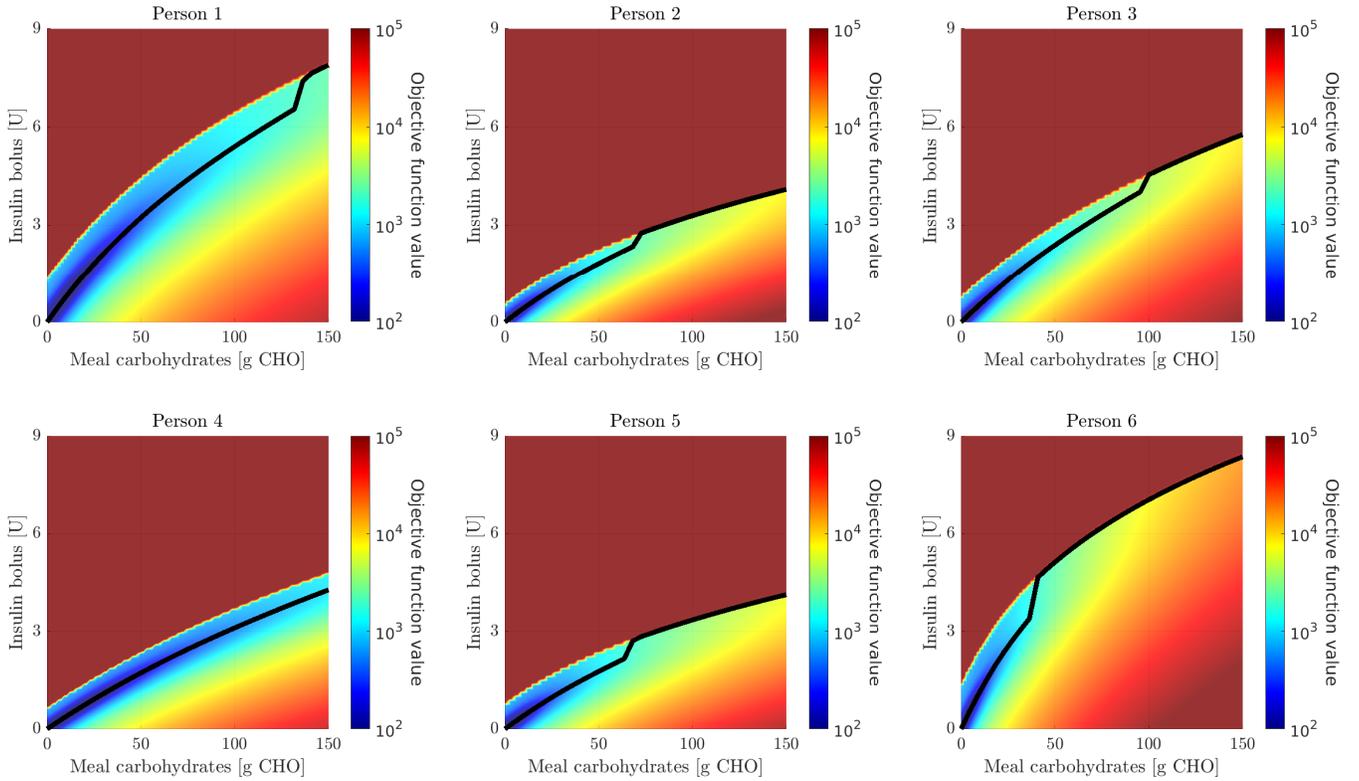}
	\caption{Values of the objective function~\eqref{eq:bolus:objective:function} for different (absolute) meal carbohydrate contents and insulin boluses. The black lines indicate the optimal insulin boluses as functions of the meal carbohydrate content, i.e., the solutions to the dynamic optimization problem~\eqref{eq:bolus}.}
	\label{fig:asymmetricQuadraticObjectiveFunctionColorPlot}
\end{figure*}

For 6~virtual people with type~1 diabetes (i.e., 6 different sets of parameters, $\theta$), Fig.~\ref{fig:asymmetricQuadraticObjectiveFunctionColorPlot} shows the values of the objective function in~\eqref{eq:bolus:objective:function} for different combinations of (absolute) meal carbohydrate contents and insulin boluses. The black lines indicate the optimal boluses found by solving the dynamic optimization problem~\eqref{eq:bolus} using a single-shooting approach. The optimal insulin bolus curve is more nonlinear for some sets of parameters than others. For instance, a linear insulin bolus curve is a worse approximation for person~6 than for person~4. The first kink in the optimal insulin bolus curve (starting from the left) appears because the total non-insulin-dependent glucose flux decreases when the blood glucose concentration comes below 4.5~mmol/L in the model by Hovorka et al.~\cite{Hovorka:etal:2004}. Consequently, more insulin is required to decrease the concentration below this value. The second kink arises because of the soft lower bound, $z_{\min}$, in the penalty function. There are no kinks for person~4 because their blood glucose concentration increases very little after meals, i.e., there would be kinks for meals with more than 150~g carbohydrates.
In conclusion, a piecewise linear function is a reasonable approximation of the optimal insulin bolus curve.
\begin{table}
	\centering
	\caption{The five glycemic ranges described by Holt et al.~\cite{Holt:etal:2021b}.}
	\label{tab:glycemic:ranges}
	\begin{tabular}{lcc}
		\hline
		Category 				& Range [mmol/L]			& Color		\\
		\hline
		Level 2 hyperglycemia 	& ]13.9, $\ \infty \; [$ 	& Orange 	\\
		Level 1 hyperglycemia 	& ]10.0, 13.9] 				& Yellow 	\\
		Normoglycemia 			& [ 3.9, 10.0] 				& Green 	\\
		Level 1 hypoglycemia 	& [ 3.0, \phantom{1}3.9[ 	& Light red \\
		Level 2 hypoglycemia 	& [ 0.0, \phantom{1}3.0[ 	& Red 		\\
		\hline
	\end{tabular}
\end{table}
	\section{Algorithm}
\label{sec:algorithm}
At time $t_k$ [min], the AP receives a CGM measurement, $y_k$ [mmol/L], and computes the basal and bolus insulin flow rates, $u_{ba, k}$ [mU/min] and $u_{bo, k}$ [mU/min]. These are clipped and collected in the vector of manipulated inputs,
\begin{align}
	u_k &=
	\max\left\{0, \min\left\{u_{\max},
	\begin{bmatrix}
		u_{ba, k} \\
		u_{bo, k}
	\end{bmatrix}
	\right\}\right\},
%
%
%
%
\end{align}
which is administered over the following control interval, $[t_k, t_{k+1}[$. The upper bounds, $u_{\max}$, on the insulin basal and bolus insulin flow rates are 55~mU/min and 8000~mU/min, respectively, and the control and sampling intervals are identical. Their lengths are $T_s = t_{k+1} - t_k = 5$~min.

If the CGM measurement is above the target value of $\bar y = 6$~mmol/L, the basal insulin flow rate is the sum of an estimated nominal basal insulin flow rate, $\bar u_{ba, k}$ [mU/min], and a microadjustment term, $u_{ma, k}$ [mU/min]. If the measurement is between the target and the safety threshold $y_s = 3$~mmol/L, the microadjustments are clipped to be non-positive (for safety reasons). Finally, if the measurement is below the safety threshold, the basal rate is zero, i.e.,
%
%
\begin{align}
	u_{ba, k} &=
	\begin{cases}
		\bar u_{ba, k} + u_{ma, k} & \text{if}~\bar y \leq y_k, \\
		\bar u_{ba, k} + [u_{ma, k}]^- & \text{if}~\bar y > y_k > y_s, \\
		0 & \text{otherwise},
	\end{cases}
\end{align}
where $[\: \cdot\: ]^- = \min\{0, \:\cdot\:\}$. The nominal basal insulin flow rate is estimated using an extended integral (I) controller, and the microadjustments are computed using a proportional-derivative (PD) controller:
\begin{subequations}
	\begin{align}
		\bar u_{ba, k} &= I_{ba, k}, \\
		u_{ma, k} &= P_{ma, k} + D_{ma, k}.
	\end{align}
\end{subequations}
Here, $I_{ba, k}$ is an integral term (see Section~\ref{sec:algorithm:estimate:basal}) and $P_{ma, k}$ and $D_{ma, k}$ are proportional and derivative terms (see Section~\ref{sec:algorithm:basal:microadjustments}).

Based on the analysis in Section~\ref{sec:analysis}, we compute the meal bolus insulin flow rate as a continuous piecewise linear function of the estimated normalized meal carbohydrate flow rate, $\hat d_k$ [g CHO/(kg min)]. That is, $\hat d_k$ is the amount of meal carbohydrates the user announces they will consume in the $k$'th control interval, divided by the product of their bodyweight and the length of the control interval, $T_s$. Specifically, the meal bolus insulin flow rate is given by
\begin{align}
	u_{bo, k} &=
	\begin{cases}
		\alpha_k d_{th} + \frac{\alpha_k}{\beta}(\hat d_k - d_{th}) & \text{if}~\hat d_k > d_{th}, \\
		\alpha_k \hat d_k & \text{otherwise}.
	\end{cases}
\end{align}
If the normalized meal carbohydrate flow rate is below the threshold $d_{th} = 0.1$~g CHO/(kg min), the insulin bolus flow rate is proportional to the meal carbohydrate flow rate and the slope, $\alpha_k$ [mU kg/(g CHO)], is essentially the inverse of the ICR. For higher normalized meal carbohydrate contents, the slope is divided by $\beta = 2$ (unitless). Both the threshold and $\beta$ were identified by trial-and-error, and we estimate the meal bolus factor $\alpha_k$ using another extended I-controller, i.e.,
\begin{align}
	\alpha_k &= I_{bo, k},
\end{align}
where $I_{bo, k}$ is an integral term described in Section~\ref{sec:algorithm:estimate:bolus}.

\subsection{Estimation of the basal rate}
\label{sec:algorithm:estimate:basal}
At time $t_k$, when a CGM measurement becomes available, we update the estimate of the basal rate and ensure that it is non-negative:
\begin{align}
	I_{ba, k} &= \max\{0, I_{ba, k-1} + \Delta I_{ba, k}\}.
\end{align}
The increment is
\begin{align}
	\Delta I_{ba, k} &= w_{ba, k} K_{I, ba} e_{ba, k} T_s,
\end{align}
where the unitless binary weight $w_{ba, k}$ is 1 if the last meal was announced more than $\Delta t_m = 9.5$~h ago and 0 otherwise. Furthermore, the integrator gain is $K_{I, ba} = 4\cdot10^{-4}$~mU L/(mmol min$^2$), and the error is computed using the deadband~$[y_{ba}^\ell, y_{ba}^u] = [3.9, 8.0]$~mmol/L and a (unitless) hypoglycemia amplification factor, $\gamma = 100$:
%
\begin{align}
	e_{ba, k} &=
	\begin{cases}
		        y_k - y_{ba}^u     & \text{if}~y_k > y_{ba}^u, \\
		\gamma (y_k - y_{ba}^\ell) & \text{if}~y_k < y_{ba}^\ell, \\
		0 & \text{otherwise}.
	\end{cases}
\end{align}

\subsection{Microadjustments of the basal rate}
\label{sec:algorithm:basal:microadjustments}
The proportional term in the microadjustment of the basal insulin flow rate is
\begin{align}
	P_{ma, k} &= w_{ma, k} K_{P, ma} e_k,
\end{align}
where the gain is $K_{P, ma} = 0.3$~mU L/(mmol min), and the error is
\begin{align}
	e_k &= y_k - \bar y.
\end{align}
%
The unitless binary weight $w_{ma, k}$ is 1 if the CGM measurement is below the target or if $w_{ba, k} = 1$. Otherwise, it is zero:
\begin{align}
	w_{ma, k} &=
	\begin{cases}
		1 & \text{if}~y_k < \bar y ~\text{or}~ w_{ba, k} = 1, \\
		0 & \text{otherwise}.
	\end{cases}
\end{align}
The derivative term is
\begin{align}
	D_{ma, k} &= w_{ma, k} K_{D, ma} \frac{y_k - y_{k-1}}{T_s},
\end{align}
where the gain is $K_{D, ma} = 10$~mU L/mmol and we disregard changes in the setpoint (which is also constant in this work).

\subsection{Estimation of the meal bolus factor}
\label{sec:algorithm:estimate:bolus}
As for the basal rate, we update the estimate of the meal bolus factor whenever a CGM measurement becomes available and ensure that it is non-negative, i.e.,
\begin{align}
	I_{bo, k} &= \max\{0, I_{bo, k-1} + \Delta I_{bo, k}\},
\end{align}
where the increment is
\begin{align}
	\Delta I_{bo, k} &= w_{bo, k} K_{I, bo} e_{bo, k} T_s.
\end{align}
The unitless binary weight $w_{bo, k}$ is 1 for a time period of $\Delta t_m$ after every announced meal, i.e., $w_{bo, k}$ and $w_{ba, k}$ never have the same value. Furthermore, the gain is $K_{I, bo} = 0.05$~mU kg L/(g CHO mmol min), and we use both a deadband of $[y_{bo}^\ell, y_{bo}^u] = [3.9, 10]$~mmol/L, the hypoglycemia amplification factor $\gamma$, and clipping to compute the error:
\begin{align}
	e_{bo, k} &=
	\begin{cases}
		        y_{bo}^{th} - y_{bo}^u  	& \text{if}~y_k > y_{bo}^{th}, \\
		        y_k         - y_{bo}^u     	& \text{if}~y_k \in [y_{bo}^u, y_{bo}^{th}], \\
		\gamma (y_k         - y_{bo}^\ell) 	& \text{if}~y_k < y_{bo}^\ell, \\
		0 & \text{otherwise}.
	\end{cases}
\end{align}
The clipping ensures that all CGM measurements above the threshold $y_{bo}^{th} = 13.9$~mmol/L result in the same error.
	\section{Virtual clinical trial}
\label{sec:numerical:example}
\begin{table}
	\setlength{\tabcolsep}{3pt}
	\caption{The compositions of the seasons and the weeks and the meal carbohydrate contents in the protocol described in~\cite{Reenberg:etal:2022}.}
	\label{tab:WeeksAndSeasons}
	\begin{tabular}{p{0.14\linewidth}|ccc}
		\multicolumn{3}{l}{\textbf{Compositions of the seasons}} \\[2pt]
		\hline
		Season & Standard weeks & Active weeks & Vacation weeks \\
		\hline
		Winter & 6 & 4 & 3 \\
		Spring & 6 & 6 & 1 \\
		Summer & 7 & 3 & 3 \\
		Autumn & 9 & 3 & 1 \\ \hline
	\end{tabular}

	\vspace{7pt}

	\begin{tabular}{p{0.14\linewidth}|cccc}
		\multicolumn{2}{l}{\textbf{Compositions of the weeks}} \\[2pt]
		\hline
		Week type & Standard days & Active days & Movie nights & Late nights \\
		\hline
		Standard & 4 & 1 & 1 & 1  \\
		Active   & 3 & 3 & 1 & 0  \\
		Vacation & 5 & 0 & 0 & 2  \\ \hline
	\end{tabular}

	\vspace{7pt}

	\begin{tabular}{l|cc}
		\multicolumn{3}{l}{\textbf{Bodyweight-dependent meal carbohydrate contents}} \\[2pt]
		\hline
		Meal size 	& Amount of carbohydrates & For a 70~kg person \\
		\hline
		Large meal  & 1.29~g CHO/kg & 90 g CHO \\
		Medium meal & 0.86~g CHO/kg & 60 g CHO \\
		Small meal  & 0.57~g CHO/kg & 40 g CHO \\
		Snack       & 0.29~g CHO/kg & 20 g CHO \\
		\hline
	\end{tabular}
\end{table}
\begin{figure}
	\centering
	\includegraphics[width=\linewidth, trim=0pt 0pt 20pt 5pt, clip]{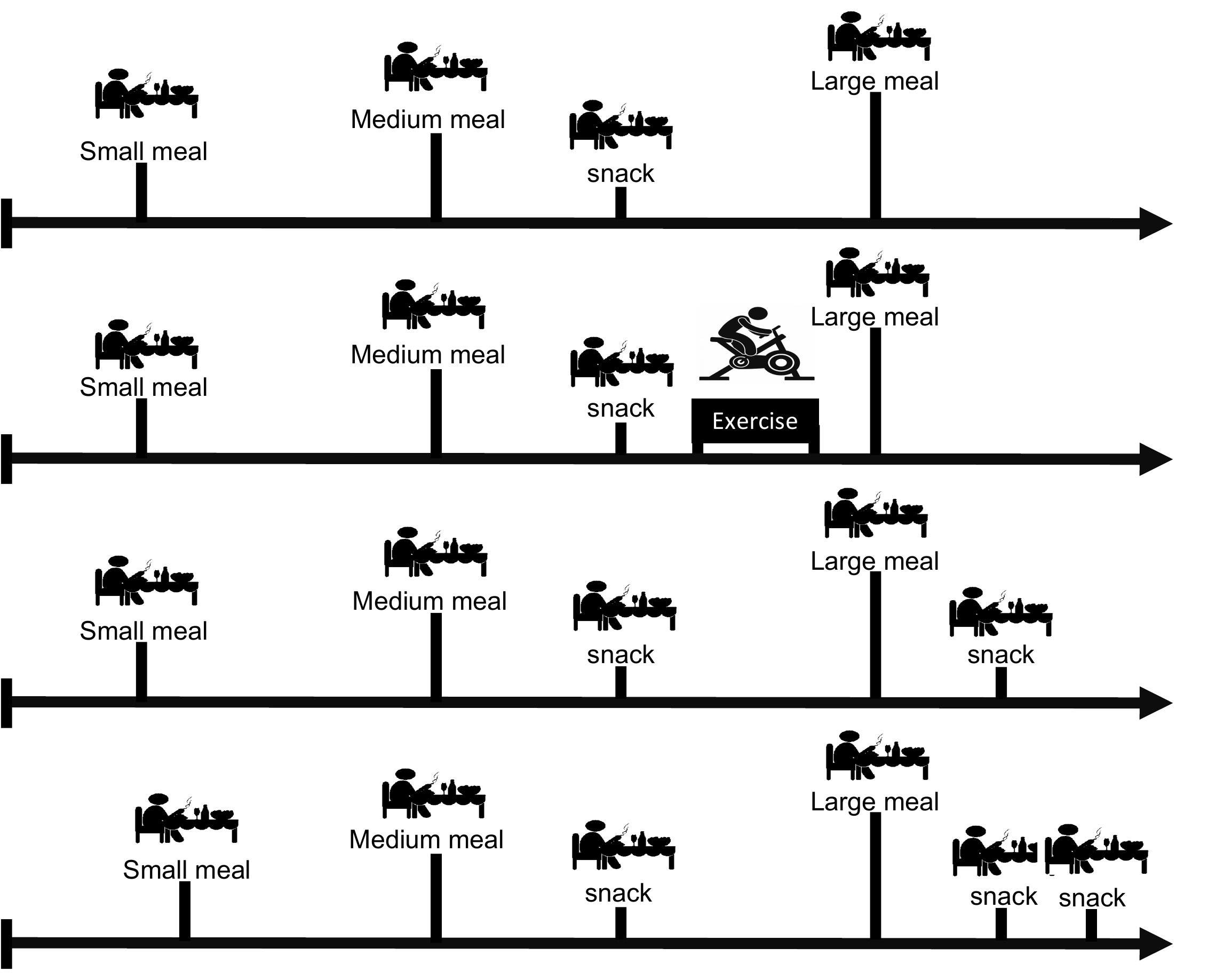}
	\caption{The different types of days in the autumn and winter of the protocol proposed by Reenberg et al.~\cite{Reenberg:etal:2022}: standard (top), active (second from the top), day with a movie night (third from the top), and day with a late night (bottom). During spring and summer, the dinner is a medium meal and the afternoon snack is consumed between breakfast and lunch instead.}
	\label{fig:WinterFallDays}
\end{figure}
\begin{figure}[tb]
	\centering
	\includegraphics[width=\linewidth, trim=5pt 75pt 85pt 64pt, clip]{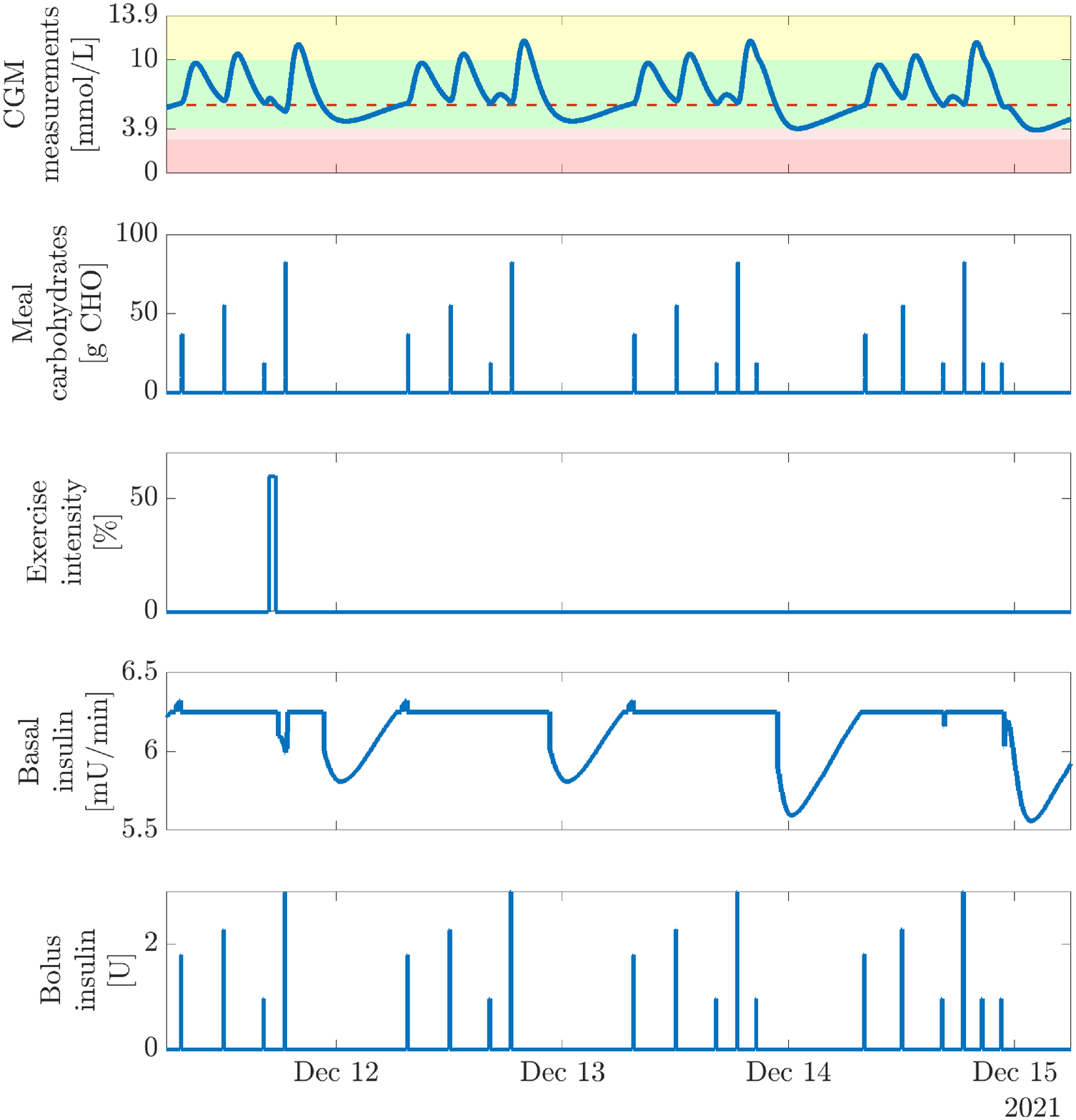}
	\caption{A single participant's CGM values (top), meal carbohydrate contents (second from the top), exercise intensity (third from the top), basal insulin flow rate (fourth from the top), and insulin boluses (bottom) over 4 days (one of each type) starting at 6:00 AM on December 11\textsuperscript{th}.
	The colored ranges are described in Table~\ref{tab:glycemic:ranges}.}
	\label{fig:CCTA2022_SinglePatientf5}
\end{figure}
In this section, we test the AP algorithm described in Section~\ref{sec:algorithm} in a virtual clinical trial with 1,000,000 participants. The trial starts on January 1\textsuperscript{st}, 2021 and lasts 52 weeks. We use 1)~the virtual population and the protocol presented by Reenberg et al.~\cite{Reenberg:etal:2022} and 2)~a previously developed Monte Carlo simulation framework~\cite{Wahlgreen:etal:2021, DTU_DCC_resource}. However, we replace participants for which any time constant is more than 1 order of magnitude smaller or larger than the mean (i.e., we generate new participants).
The protocol mimics a Northern European lifestyle, and it consists of 4 seasons lasting 13 weeks each. Each week is categorized as \emph{standard}, \emph{active}, or \emph{vacation}, and all weeks consist of \emph{standard} days, \emph{active} days, days with a \emph{movie night}, and days with a \emph{late night} (see Table~\ref{tab:WeeksAndSeasons} and Fig.~\ref{fig:WinterFallDays}). Each participant is represented using the mathematical model presented by Hovorka et al.~\cite{Hovorka:etal:2004} extended with a CGM model~\cite{Facchinetti:etal:2014} and an exercise model~\cite{Rashid:etal:2019}. The initial estimates of the basal rate and the meal bolus factor are zero, i.e., $I_{ba, 0} = I_{bo, 0} = 0$, and the initial state is the steady state without insulin administration.

Fig.~\ref{fig:CCTA2022_SinglePatientf5} shows the results of the virtual clinical trial for one participant over four different types of days. The basal rate is constant for most parts of the day, and it is decreased when the CGM measurements are below the target of 6~mmol/L (which mostly happens at night). A bolus is administered for each meal, and for this participant, the majority of the CGM measurements are within the normoglycemic range. However, the estimated nominal basal rate is quite low. Therefore, the CGM measurements increase over night.

\begin{figure}
	\includegraphics[width=0.4\linewidth, trim=0pt 50pt 35pt 35pt, clip]{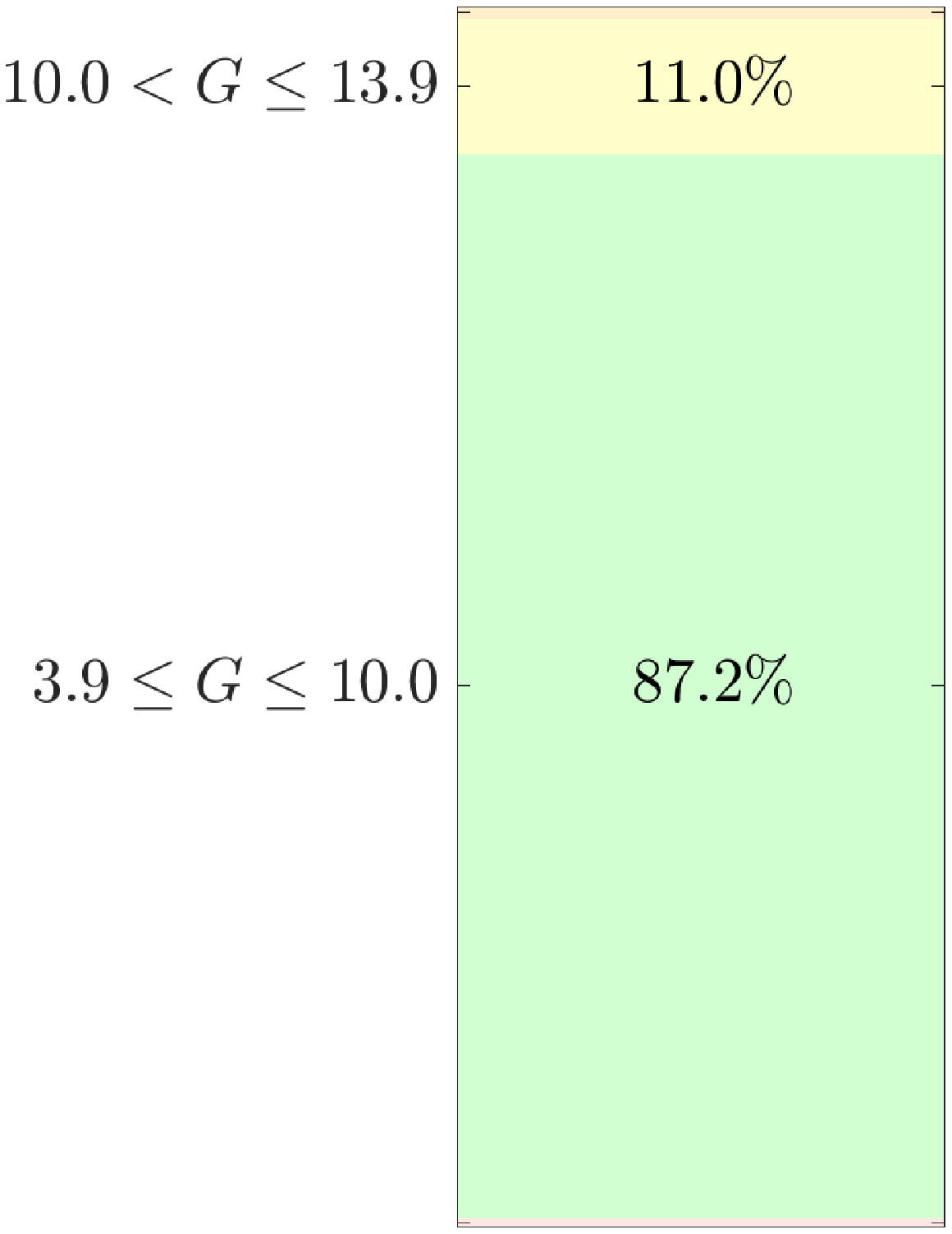}~
	\includegraphics[width=0.4\linewidth, trim=0pt 50pt 35pt 35pt, clip]{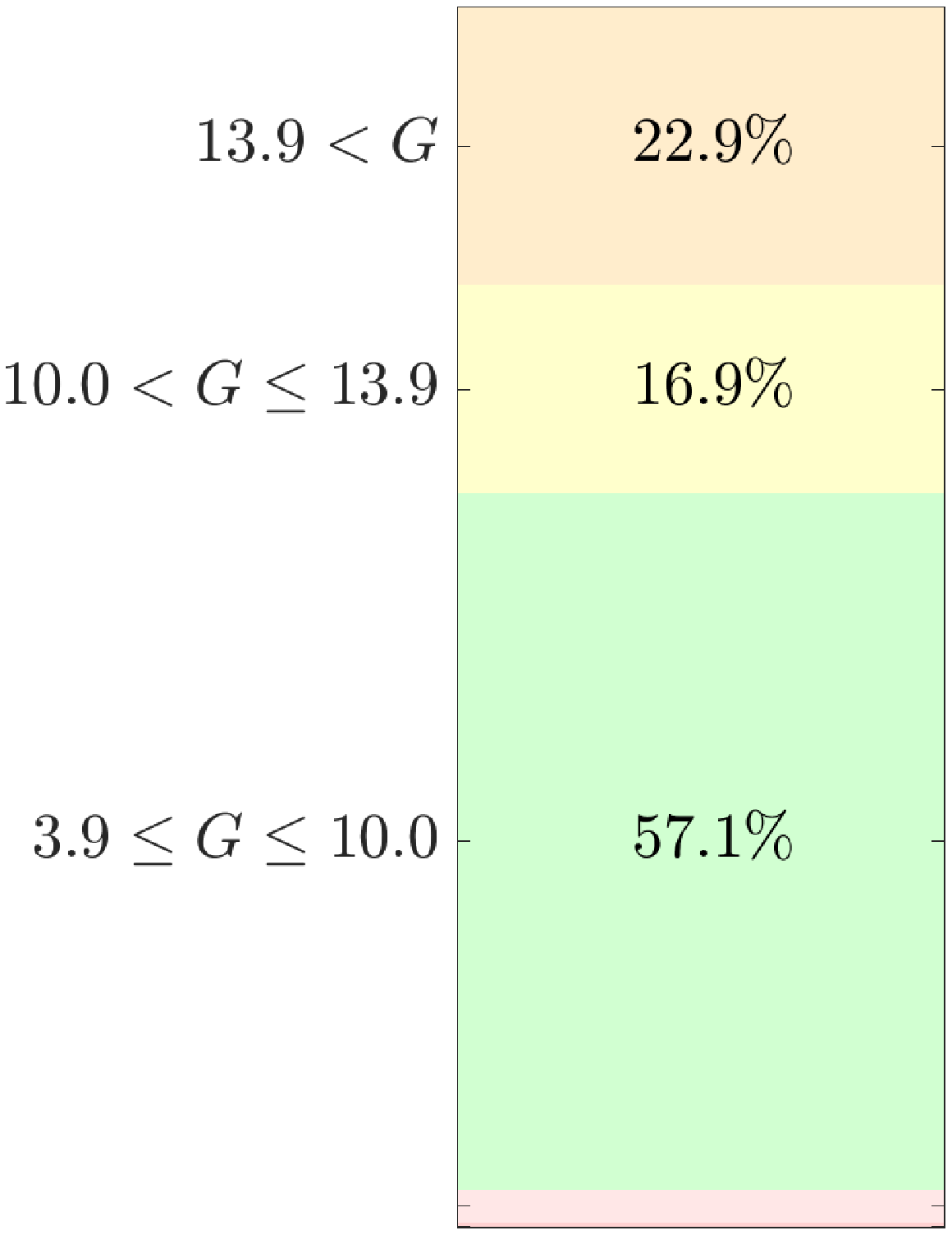}
	\caption{The mean TIRs (left) and the TIRs for the worst-case participant (right) based on the ranges in Table~\ref{tab:glycemic:ranges}.}
	\label{fig:StackedBarPlots}
\end{figure}

\begin{figure}
	\includegraphics[width=\linewidth, trim=45pt 0pt 75pt 25pt, clip]{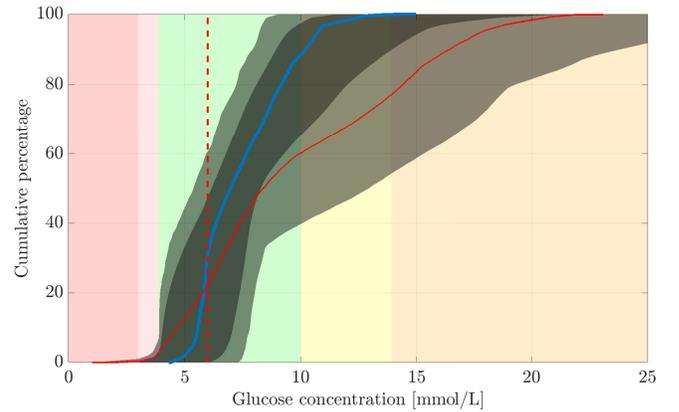}
	\caption{The cumulative distribution of the CGM measurements for the mean (blue solid line), the worst-case participant (red solid line), the 95\% central range (dark grey shaded area), and for all participants (light grey shaded area). The target is 6~mmol/L (red dashed line). The colored ranges are described in Table~\ref{tab:glycemic:ranges}.}
	\label{fig:CumulativePlot}
\end{figure}

\begin{figure}[t]
	\includegraphics[width=\linewidth, trim=75pt 30pt 100pt 40pt, clip]{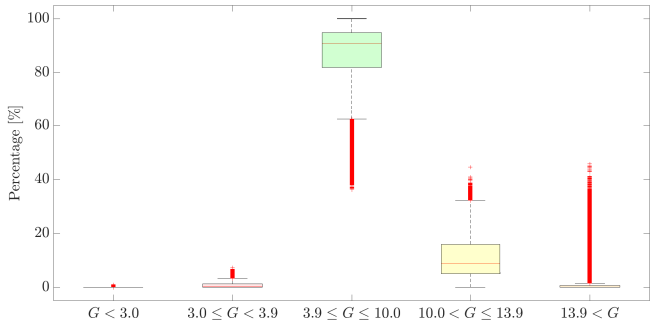}
	\caption{Box plots of the TIRs with medians (red horizontal lines), boxes spanning the first to the third quartile, and whiskers (solid black horizontal lines).
	The whiskers are 1.5 times the interquartile ranges (the height of the boxes) above or below the medians, unless the most extreme values are closer to the medians. In that case, the whiskers are the most extreme values. The red pluses are values that are beyond the whiskers (i.e., outliers).}
	\label{fig:BoxPlots}
\end{figure}

\begin{figure}
	\includegraphics[width=\linewidth, trim=60pt 10pt 85pt 25pt, clip]{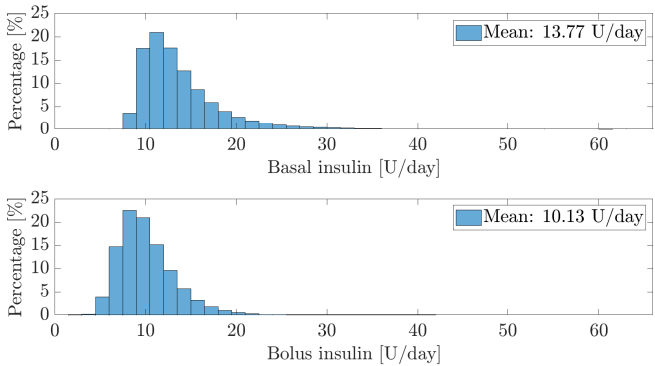}
	\caption{The distributions of the average TDDs of basal and bolus insulin. Both distributions have long right tails which are hardly visible.}
	\label{fig:TotalDailyDose}
\end{figure}

In the following, we discuss the efficacy of the AP based on the last 48~weeks of the trial as the estimates of the nominal basal insulin flow rate and the meal bolus factor vary significantly during the first 4 weeks.
The participant who obtains the lowest CGM measurement during all 52 weeks (specifically, 1.05~mmol/L) is referred to as the worst-case participant, and Fig.~\ref{fig:StackedBarPlots} shows the mean and worst-case TIRs.  The mean of 87.2\% TIR exceeds the target of 70\%, and the time in level 1 and 2 hypoglycemia is low, even for the worst-case participant.
Fig.~\ref{fig:CumulativePlot} shows the cumulative distribution of the CGM measurements. The left tail shows that all participants spend less than 1\% of the time in level 2 hypoglycemia and less than 8\% of the time in level 1 and 2 hypoglycemia.
This can also be seen in Fig.~\ref{fig:BoxPlots} which shows box plots of the TIRs. It also shows that most participants do not spend significant amounts of time in level 2 hyperglycemia.
Table~\ref{tab:diabetes:treatment:targets} shows the percentages of participants satisfying the glycemic targets described by Holt et al.~\cite{Holt:etal:2021b}. Almost 82\% satisfy all targets, and nearly 89\% satisfy all average glucose, GMI, TAR, TIR, and TBR targets.
Finally, Fig.~\ref{fig:TotalDailyDose} shows that most of the participants' average total daily doses (TDDs) of basal and bolus insulin are in the intervals [7.5, 25]~U/day and [5, 20]~U/day, respectively. However, the distributions have long tails towards the right indicating that a few participants require high insulin doses.
\begin{table}[t]
	\centering
	\caption{Satisfaction of the glycemic targets described in~\cite{Holt:etal:2021b}.}
	\label{tab:diabetes:treatment:targets}
	\begin{tabular}{llr}
		Quantity 									& Target 		& Satisfied \\
		\hline \\[-7pt]
		Average glucose 							& $< 154$~mg/dL &  97.69\% 	\\
		GMI 										& $< 7\%$ 		&  97.77\% 	\\
		GV 											& $\leq 36\%$ 	&  82.70\% 	\\
		\hline \\[-7pt]
		TAR (level 2 hyperglycemia) 				& $< 5\%$ 		&  93.97\% 	\\
		TAR (level 1 and 2 hyperglycemia) 			& $< 25\%$ 		&  89.06\%  \\
		TIR\; (normoglycemia)						& $> 70\%$ 		&  92.19\%  \\
		TBR (level 1 and 2 hypoglycemia) 			& $< 4\%$ 		&  98.85\%  \\
		TBR (level 2 hypoglycemia) 					& $< 1\%$ 		& 100.00\%  \\
		\hline \\[-7pt]
		\textbf{All TAR, TIR, and TBR targets} 		& 				& 88.92\% 	\\
		\textbf{All targets except the GV target} 	& 				& 88.87\% 	\\
		\textbf{All targets} 						& 				& 81.70\% 	\\
		\hline
	\end{tabular}
\end{table}
	\section{Conclusions}
\label{sec:conclusions}
We propose a one-size-fits-all AP algorithm for people with T1D, which estimates both the basal insulin flow rate and the meal insulin bolus curve. It is based on physiological insight and concepts from PID control, and it only requires the bodyweight, CGM measurements, and meal carbohydrate estimates.
We compute the meal insulin bolus as a piecewise linear function of the meal carbohydrate content normalized with bodyweight, and we test the AP algorithm in a 52 week virtual clinical trial with 1,000,000 participants. The mean TIR is 87.2\%, and almost 89\% of the participants satisfy targets on average glucose, GMI, TAR, TIR, and TBR.

	\bibliographystyle{IEEEtran}
	\bibliography{./ref/References}
\end{document}